\documentclass[12pt]{amsart}
\usepackage{amsmath,amsthm,amsfonts,amssymb}
\begin{document}
\hyphenation{Ko-ba-ya-shi}
\newcommand{\Aut}{\mathop{\rm Aut}}
\newcommand{\A}{{\mathbb A}}
\newcommand{\B}{{\mathbb B}}
\newcommand{\C}{{\mathbb C}}
\newcommand{\N}{{\mathbb N}}
\newcommand{\Q}{{\mathbb Q}}
\newcommand{\Qp}{{\mathbb Q}_+}
\newcommand{\Z}{{\mathbb Z}}
\renewcommand{\O}{{\mathcal O}}
\renewcommand{\P}{{\mathbb P}}
\newcommand{\R}{{\mathbb R}}
\newcommand{\rc}{\subset}
\newcommand{\rank}{\mathop{rank}}
\newcommand{\Sing}{\mathop{Sing}}
\newcommand{\tensor}{\otimes}
\newcommand{\mult}{\mathop{mult}}
\newcommand{\trace}{\mathop{tr}}
\newcommand{\dimc}{\mathop{dim}_{\C}}
\newcommand{\Lie}{\mathop{Lie}}
\newcommand{\Auto}{\mathop{{\rm Aut}_{\mathcal O}}}
\newcommand{\alg}[1]{{\mathbf #1}}
\newtheorem{definition}{Definition}
\newtheorem*{claim}{Claim}
\newtheorem{corollary}{Corollary}
\newtheorem{Conjecture}{Conjecture}
\newtheorem*{SpecAss}{Special Assumptions}
\newtheorem{example}{Example}
\newtheorem*{remark}{Remark}
\newtheorem*{notation}{Notation}
\newtheorem*{observation}{Observation}
\newtheorem*{fact}{Fact}
\newtheorem*{remarks}{Remarks}
\newtheorem{lemma}{Lemma}
\newtheorem{proposition}{Proposition}
\newtheorem{theorem}{Theorem}
\title{%
A Brody theorem for orbifolds
}
\author {Frederic Campana \&\ J\"org Winkelmann}
\begin{abstract} We study the Kobayashi pseudodistance
for orbifolds, proving an orbifold version of Brody's theorem
and classifying which one-dimensional orbifolds are hyperbolic.
\end{abstract}
\address{%
 Frederic Campana \& J\"org Winkelmann \\
 Institut Elie Cartan (Math\'ematiques)\\
 Universit\'e Henri Poincar\'e Nancy 1\\
 B.P. 239\\
 F-54506 Vand\oe uvre-les-Nancy Cedex\\
 France
}
%
\maketitle
\section{Introduction}
We study orbifolds as introduced in \cite{C}, define morphisms
and discuss hyperbolicity. For this purpose we establish
a Brody theorem for orbifolds (see \cite{B} for the Brody theorem
for complex spaces (compare also \cite{Z} for a different approach)).
Using this Brody theorem for orbifolds we then determine
which one-dimensional orbifolds are hyperbolic.

There are two different classes of orbifold morphisms, baptised
``classical'' resp.~``non-classical''.

In the ``classical sense''
many problems are easier to handle because ``classical''
orbifold morphisms
behave very well with respect to \'etale orbifold morphisms.
In particular, the
classification of one-dimensional hyperbolic orbifolds
can be obtained via ``unfoldings''.

In contrast, for determining which one-dimensional orbifolds
are hyperbolic in the ``non-classical'' sense we really need
our ``Brody theorem for orbifolds''.

\section{Orbifolds}
We always assume all complex spaces to be irreducible, reduced, normal,
Hausdorff and paracompact.

We recall some notions introduced in \cite{C}.

Let $\Qp=\{x\in\Q:x> 0\}$.
An {\em effective
Weil $\Qp$-divisor} on a complex space $X$ is a formal sum
$\Delta=\sum a_i[Z_i]$ with all its coefficients
$a_i$ in $\Qp$, the $Z_i$ being pairwise distinct
irreducible reduced hypersurfaces on $X$. The support $\vert \Delta\vert$ of $\Delta$ is the union of all $Z_i's$.
An orbifold $(X/\Delta)$ is a pair consisting of an
irreducible complex space $X$ together with a Weil 
$\Qp$-divisor $\Delta=\sum a_i[Z_i]$ for which 
\[
a_i\in\{1\}\cup\{1-\frac{1}{m}: m \in\N\}\ \forall i
\]

In this case, $a_i=1-1/m_i$ (resp $m_i=1/(1-a_i)$) is the weight (resp. the multiplicity) of $Z_i$ in $\Delta$.
It is convenient to consider
$\infty$ as the multiplicity for the weight $1=1-\frac{1}{\infty}$.

If $\Delta$ is an empty divisor, we will frequently identify 
$(X/\Delta)$ with $X$.
If $Z$ is a component of weight $1$ of $\Delta$, then
$(X/\Delta)$ may (and frequently will be) identified with
$(X'/\Delta')$ where $X'=X\setminus Z$ and $\Delta'=\Delta-[Z]$.

An orbifold $(X/\Delta)$ is called {\em compact} iff
$X$ is compact and $\Delta$ contains no irreducible component
of multiplicity $1$ (=weight $\infty$).

An orbifold $(X/\Delta)$ is {\em smooth} (or non-singular)
if $Y$ is smooth {\em and} $\Delta$ is a locally s.n.c.~divisor.

\section{Orbifold morphisms}
Orbifold were introduced in \cite{C} in the context of fibrations.
For a reducible fiber of a fibration there are two ways to define
its multiplicity: Classically one takes the greatest common divisor
of the multiplicities of its irreducible components. Non-classically
(and this is the point of view emphasized in \cite{C}) one takes
the infimum of these multiplicities. Correspondingly, we define
two notions of orbifold morphisms, a ``classical'' one and a
``non-classical'' one.

\begin{definition}\label{def-disc-map}
Let $(X/\Delta)$ be an orbifold with $\Delta=\sum_i (1-\frac{1}{m_i})Z_i$
where $m_i\in\N\cup\{\infty\}$ and where the $Z_i$ are distinct
irreducible hypersurfaces.

A holomorphic map $h$
from the unit disk $D=\{z\in\C:|z|<1\}$ to $X$ is
a (non-classical) {\em ``orbifold morphism from $D$ to $(X/\Delta)$''} if
$h(D)\not\subset|\Delta|$ and if, moreover
$\mult_x(h^*Z_i)\ge m_i$ for all $i$ and $x\in D$ with
$h(x)\in|Z_i|$. If $m_i=\infty$, we require $h(D)\cap Z_i={\emptyset}$.
The map $h$ is called a ``classical orbifold morphism''
if the condition ``$\mult_x(h^*Z_i)\ge m_i$'' is replaced by
the condition ``$\mult_x(h^*Z_i)$ is a multiple of $m_i$''.
\end{definition}

\begin{definition}\label{def-orbi-map}
Let $(X/\Delta)$ and $(X'/\Delta')$ be orbifolds.
Let $\Delta_1$ be the union of all irreducible components of $\Delta$
with multiplicity $1$ (equivalently: weight $\infty$).
An {\em ``orbifold morphism''}
(resp.~``classical orbifold morphism'')
from $(X/\Delta)$
to $(X'/\Delta')$ is a holomorphic map $f:X\setminus\Delta_1\to X'$
such that
\begin{enumerate}
\item
$f(X)\not\subset|\Delta'|$.
\item
For every orbifold morphism
(resp.~every classical orbifold morphism)
in the
sense of def.~\ref{def-disc-map}, 
$g:D\to (X/\Delta)$ with $g(D)\not\subset f^{-1}(|\Delta'|)$
the composed map $f\circ g:D\to X'$
defines an orbifold morphism from $D$ to $(X'/\Delta')$.
\end{enumerate}
\end{definition}

\begin{remark}
If one would like to obtain a purely algebraic-geometric definition
in the case where $X$, $X'$, $\Delta$ and $\Delta'$ are algebraic,
one might use smooth algebraic curves instead of the unit disc.
\end{remark}

\begin{remark} 
Both notions (``classical'' versus ``non-classical'')
differ substantially. The class of ``classical orbifold morphism'' is
a much more restricted one. Questions concerning ``classical'' morphism
often can be handled easily by using \'etale orbifold morphisms,
which does not work in the non-classical setup.

Similar differences occur also
for questions related to function field and arithmetic versions,
see \cite{C'}.
\end{remark}

\section{Examples of orbifold morphisms}

\subsection{Elementary properties}

{\em Composition.}
If $f:(X/\Delta)\to (X'/\Delta')$ and $g:(X'/\Delta')\to (X''/\Delta'')$
are orbifold morphisms resp.~classical orbifold morphisms,
 so is $g\circ f:(X/\Delta)\to (X''/\Delta'')$
unless $(g\circ f)(X)\subset|\Delta''|$.

{\em Empty divisors}.
If $\Delta$, $\Delta'$ are empty Weil divisors on complex
spaces $X$ resp.~$X'$, then every holomorphic map from $X$ to $X'$
defines a (classical) orbifold morphism from $(X/\Delta)$ to $(X'/\Delta')$.

{\em Majorisation/Minorisation.}
If $f:(X/\Delta)\to(X'/\Delta')$ is an orbifold morphism
and $\Delta''$ is a $\Qp$-Weil divisor on $X'$ with
$\Delta''\le\Delta'$, then $f$ is an orbifold morphism
to $(X'/\Delta'')$, too.

Similarly: If $f:(X/\Delta)\to(X'/\Delta')$ is an orbifold morphism
and $\Delta''$ is a $\Qp$-Weil divisor on $X$ with
$\Delta''\ge\Delta$, then $f$ is an orbifold morphism
from $(X/\Delta'')$, too.

If $f:X\to Y$ is a holomorphic map of complex spaces and $D$ is
an irreducible reduced hypersurface on $Y$ such that
\[
f:X\to \left( Y/ (1-\frac{1}{n})[D]\right)
\]
is an orbifold morphism for all $n\in\N$,
then $f(X)\cap|D|=\emptyset$.

\subsection{Curves}
Let $C$ and $C'$ be smooth complex curves, $p\in C$, $p'\in C'$, $n,n'\in\N$.
Then a non-constant holomorphic map $f:C\to C'$ is an orbifold morphism
from $(C/(1-\frac{1}{n})\{p\})$
to $(C'/(1-\frac{1}{n'})\{p'\})$
if  $\mult_pf^*[\{p'\}]\ge \frac{n'}{n}$
and $\mult_zf^*[\{p'\}]\ge n'$ for $z\in C\setminus\{p\}$.

In particular, if $C=C'$ and $p=p'$, then the identity map
defines an orbifold morphism iff $n\ge n'$.

In addition, $f:C\to C'$ is a ``classical orbifold morphism'' if 
``$\ge$'' is replaced by ``{\em is a multiple of}'', i.e.
$n\mult_pf^*[\{p'\}]$ and $\mult_zf^*[\{p'\}]$ for $z\ne p$
 must divide $n'$.

\subsection{Automorphisms}
Let $(X/\Delta)$ be an orbifold.
A holomorphic automorphism $f$
of $X$ is an orbifold morphism iff $f^*\Delta=\Delta$.

\subsection{Blown up surface}
Let $S$ be a complex surface and $\pi:\hat S\to S$
the sigma-process centered at a point $c\in S$.
Let $D_i$ be a finite family of irreducible reduced
hypersurfaces (i.e.~curves) on $S$ with total
transforms $\pi^*D_i$ and strict transforms $\hat D_i$.
Then $\pi^*D_i=\hat D_i+d_iE$ where $E=\pi^{-1}(c)$ and
where $d_i$ denotes the multiplicity of $D_i$ at $c$.
Let $\Delta=\sum_i (1-\frac{1}{n_i})D_i$ for some $n_i\in\N$
and let $\hat\Delta$ be a $\Qp$-Weil divisor on $\hat S$.
Let $m=\max_i\frac{n_i}{d_i}$.
Then $\pi:\hat S\to S$ defines an orbifold morphism
from $(\hat S/\hat\Delta)$ to $(S/\Delta)$ iff
each $\hat D_i$ occurs with multiplicity at least $(1-\frac{1}{n_i})$
in $\hat\Delta$ and in addition $E$ occurs with multiplicity at least
$ 1-\frac{1}{m}$ in $\hat\Delta$.

\subsection{Quotients by group actions}
Let $G$ be a discrete group acting effectively on a complex curve $Y$.
Such an action is called ``proper'' resp.~``properly discontinuously''
if the map $\mu:G\times Y\to Y\times Y$ given by
$\mu(g,y)=(g\cdot y,y)$ is a proper map.
In particular, if $G$ is finite, then every action of $G$ is proper.
The quotient $X=Y/G$ has a the structure of a ringed topological
space in a canonical way. If $G$ is acting properly and $Y$ is smooth, 
then $Y/G$ is a smooth complex curve.

For $y\in Y$ let $G_y$ denote the isotropy group at $y$,
i.e.~$G_y=\{g:g\cdot y=y\}$.
Assume that $\dim(Y)=1$. In this case $Y/G$ is smooth and furthermore
we can define a $\Qp$-divisor $\Delta$ on $X=Y/G$
by $\Delta=\sum_{[y]\in Y/G}(1-1/\#G_y)\{[y]\}$.

Then $(X/\Delta)$ is an orbifold such that the natural projection
from $Y$ onto $(X/\Delta)$ is an orbifold morphism.

Moreover this orbifold morphism is \'etale in the sense of
definition~\ref{def-etale}.

\section{Ramification divisors}

\subsection{Existence}

\begin{theorem}\label{ex-ram}
Let $f:X\to Y$ be a surjective holomorphic map with constant fiber dimension
between irreducible normal complex spaces.

Then there exists a
unique Weil divisor $R_f$ on $X$ with the following properties:
\begin{enumerate}
\item If $D_0$ and $D_1$ are reduced irreducible hypersurfaces on $X$
resp.~$Y$ and $D_0$ occurs with multplicity $m\ge 2$ in $f^*D_1$,
then $D_0$ occurs with multiplicity $(m-1)$ in $R_f$.
\item $f(R_f)$ contains no open subset of $Y$.
\end{enumerate}
\end{theorem}

\begin{notation}
This divisor $R_f$ is called {``ramification divisor''}.
\end{notation}
\begin{remark}
If the map is not surjective or the fibers are not equidimensional,
then in general there is no such divisor with these properties.
\end{remark}
\begin{proof}
We simply define $R_f$ as the sum of all $D_0$ with respective
multiplicities as required by the first property.
There are two problems in doing so:
\begin{itemize}
\item
Given an irreducible reduced hypersurface $D_0\subset X$,
we need that there is at most one irreducible reduced
hypersurface $D_1\subset Y$ such that $|D_0|\subset|f^*D_1|$.
\item
The sum must be locally finite.
\end{itemize}
The first property is a consequence of the assumption that $f$ is
surjective with equidimensional fibers.
For the second we observe that, for any such $D_0$ with multiplicity $\ge 2$,
the support $|D_0|$ must not intersect the set
$\Omega$ of all non-singular points $x\in X$ for which $f(x)$
is non-singular and $Df:T_xX\to T_{f(x)}Y$ is surjective.
The complement of $\Omega$ is an analytic subset of $X$, hence
it locally contains only finitely many hypersurfaces.
For this reason the sum of all such $D_0$ is locally finite.
\end{proof}

\begin{proposition}
If $f$ is a surjective finite morphism between complex manifolds
$X$ and $Y$, then $R_f$ is linearly equivalent to $K_X\tensor (f^*K_Y)^*$.
\end{proposition}
This follows by pulling-back $n$-forms ($n=\dim(X)=\dim(Y)$).

There is no such statement in the case where the fibers are
positive-dimensional: Let $C$ be a compact smooth curve
and let $p_1$, $p_2$ be the projections from the
product $X=C\times\P_1$ to its factors.
Then $R_{p_1}=0$, but $K_X\tensor(p_1^*K_C)^{-1}\sim p_2^*K_{\P_1}$.

\subsection{ Composition rule.}

\begin{proposition}
Let $f:X\to Y$ and $g:Y\to Z$ be surjective holomorphic maps
with equidimensional fibers between normal complex spaces.
Then
\[
R_{g\circ f}=R_f+f^*R_g-S_{f,g}
\]
where
$S_{f,g}$ denotes the sum of those irreducible components
of $R_f$ which are mapped dominantly on $Y$ by $g\circ f$.
\end{proposition}

\subsection{Orbifold morphisms and ramification divisor.}

\begin{proposition}
Let $(X/\Delta)$ and $(Y/\Delta')$ be smooth orbifolds.
Let $f:X\to Y$ be a surjective holomorphic map with constant fibre dimension
between irreducible complex spaces.

Then $f$ defines an orbifold morphism from $(X/\Delta)$
to $(Y/\Delta')$ if and only if $(R_f+\Delta-f^*\Delta')\ge 0$.
\end{proposition}
\begin{proof}
We may check this for each irreducible component separately.
Thus let $H$ be an irreducible component of $\Delta$ with
multiplicity $(1-1/n)$ and let $H'$ be an irreducible component
of $\Delta'$ with multiplicity $(1-1/m)$ such that
$|H|\subset|f^*H'|$. Assume that $H$ occurs with multiplicity
$d$ in $f^*H'$.

In order for to be an orbifold morphism, we need that
$g^*f^*H'$ has multiplicity at least $m$ whenever $g:D\to X$
is a holomorphic map for which $g^*H$ has multiplicity $\ge n$.
This is the case if $nd\ge m$.

On the other hand the multiplicity of $H$ in $R_f+\Delta-f^*\Delta'$
equals $(d-1)+(1-1/n)-(d(1-1/m))$.
Now
\[
(d-1)+(1-1/n)-(d(1-1/m))= -1/n + d/m
\]
and $(-1/n+d/m)\ge 0$ holds if and only if $nd\ge m$.
\end{proof}

\begin{lemma}
Assume that there exists a
non-constant orbifold morphism $f:(C/\Delta)\to(C'/\Delta')$
for some smooth compact Riemann surfaces $C$ and $C'$.
Let $K_C$ and $K_{C'}$ denote the respective canonical line bundles
on $C$ resp.~$C'$.

Then
\[
\deg(K_C+\Delta) \ge d.\deg(K_{C'}+\Delta'),
\] if $d$ is the geometric degree of $f$ (ie: the number of points of one of its generic fibres).
\end{lemma}
\begin{proof}
Because $f$ is an orbifold morphism, we have $R_f+\Delta-f^*\Delta'\ge 0$.
On the other hand, $R_f\sim K_C-f^*K_{C'}$.
Therefore
\[
\deg\left(K_C-f^*K_{C'}+\Delta-f^*\Delta'\right)\ge 0.
\]
Hence
\[
\deg\left(K_C+\Delta\right) \ge
(\deg f)\deg\left(K_{C'}+\Delta'\right)\ge
\deg\left(K_{C'}+\Delta'\right).
\]
\end{proof}

\section{Orbifold base}

\begin{lemma}\label{max}
Let $(X/\Delta)$, $(Y/\Delta')$ and $(Y/\Delta'')$
be orbifolds.

Assume that $f:X\to Y$ is a holomorphic map which defines an orbifold
morphism from $(X/\Delta)$ to both $(Y/\Delta')$ and $(Y/\Delta'')$.

Then $f$ likewise defines an orbifold morphism to
$(Y/\max\{\Delta',\Delta''\})$.
\end{lemma}

\begin{proof}
Immediate.
\end{proof}

\begin{definition}\label{def-orbifold-base}
Let $f:X\to Y$ be a holomorphic map of complex spaces.
Then $(Y/\Delta)$ is an ``orbifold base'' for $f$ if
$\Delta$ is a maximal $\Qp$-Weil divisor for which
$f$ defines an orbifold morphism from $(X/\emptyset)$ to $(Y/\Delta)$.
\end{definition}

In view of lemma~\ref{max}
the following is immediate:

\begin{lemma}
Let $f:X\to Y$ be a holomorphic map of complex spaces.

Either there exists an orbifold base or there is an infinite
sequence of distinct irreducible reduced hypersurfaces $H_i$ on $Y$
such that $f:X\to (Y/\frac{1}{2}H_i)$ is an orbifold morphism.
\end{lemma}

\begin{proposition}
There exists an orbifold base $(Y/\Delta)$ for every surjective holomorphic
map $f$ between irreducible reduced complex spaces $X$, $Y$.
\end{proposition}

\begin{proof}
Let $H$ be an irreducible reduced hypersurface in $Y$ for which
there exists a number $n\ge 2$ such that
$f$ is an orbifold morphism to $(Y/(1-\frac{1}{n}))H$.
Then for every $p\in X$, $q=f(p)\in H$ and every holomorphic
map $g:D\to X$ with $$g(0)=p$$ we have $\mult_0((f\circ g)^*D\ge n\ge 2$.

Let \[
\Omega=\{x\in X_{reg}:Tf_x \text{ is surjective }\}.
\]
Then $\Omega$ can not intersect $f^{-1}(H_{reg})$.
Therefore $|H|\subset Y\setminus f(\Omega)$.
But $Y\setminus f(\Omega)$ is an analytic subset of $Y$.
It follows that the family of all hypersurfaces $H_i$
for which there exists a number $n_i$ such that
$f:X \to(Y/(1-\frac{1}{n_i})H_i)$ is a locally finite family.
Hence
\[
\Delta =\max (1-\frac{1}{n_i})H_i
\]
exists and $(Y/\Delta)$ is the orbifold base for $f:X\to Y$.
\end{proof}

\begin{remark}
Surjectivity of $f$ is crucial, as shown by the following example
of a curve $Q$ and a holomorphic map $i:Q\to\P_2$ for which
there are infinitely many curves $L_s$ in $\P_2$ such that
$f:Q\to (\P_2/\frac{1}{2}L_s)$
is an orbifold morphism.

Let $S$ be a finite subset of a smooth quadric $Q$ in $\P_2$.
For each $s\in S$ let $L_s$ denote the line through $s$ which
is tangent to $Q$ at $s$. Since $\deg(Q)=2$, the two curves $Q$and $L_s$ intersect only at $s$ and there with multiplicity two.
Then the embedding $i:Q\to \P_2$ defines an orbifold
morphism from $Q=(Q/\emptyset)$ to $(\P_2/\Delta)$
with\[\Delta=\sum_{s\in S}\frac{1}{2}[L_s]\]
Note that $S$ is an arbitrary finite subset, we do not need any bound on its cardinality.

\end{remark}

\section{Canonical divisors}

\begin{definition}
For a smooth orbifold $(X/\Delta)$ we define the canonical divisor
$K_{(X/\Delta)}$ as
$K_X+\Delta$.

$(X/\Delta)$ is said to be of ``general type'' if
$K_{(X/\Delta)}=K_X+\Delta$ is a big divisor on $X$.
\end{definition}

(A $\Q$-divisor $D$ is called ``big'' if there exists a natural number 
$n$ such that $nD$ is a ($\Z$-)divisor and the sections of 
the associated line bundle $L(nD)$ yield a bimeromorphic map
from $X$ to a subvariety of $\P(\Gamma(X,L(nD))^*)$.)

\subsection{Etale morphisms}

\begin{definition}\label{def-etale}
Let $(X/\Delta)$ and $(X'\Delta')$ be smooth orbifolds.
An orbifold morphism $\pi:(X/\Delta)\to (X'\Delta')$ is called
{\em \'etale} if the following two conditions are fulfilled:
\begin{enumerate}
\item
the fibers of the underlying map $\pi:X\to X'$ are discrete,
\item
the underlying map $\pi:X\to X'$ is weakly proper, 
i.e., for every point $p\in X'$ 
there is an open neighbourhood $U(p)$ such that the restriction of $\pi$
to any connected component of $\pi^{-1}(U)$ is proper in the usual sense,
\item and $R_\pi=\Delta-\pi^*\Delta'$ where $R_\pi$ is the ramification
divisor of $\pi:X\to X'$.
\end{enumerate}
\end{definition}

Note that $R_\pi$ exists by thm.~\ref{ex-ram} in view of the first condition.

If $(X/\Delta)$ and $(X'/\Delta')$ are compact, this is equivalent
to the condition $\pi^*K_{(X'/\Delta')}\simeq K_{(X/\Delta)}$.

If in addition $X$ and $X'$ are one-dimensional,
a finite morphism $f:X\to X'$ defines an \'etale orbifold
morphism
if and only if
it is an orbifold morphism, and:
\[
\deg(K_X+\Delta) = \deg(K_{(X/\Delta)})=
d.\deg(K_{X'}+\Delta')=d.\deg(K_{(X'/\Delta')}).
\]

A holomorphic map $f:D\to D$ gives an \'etale orbifold morphism
from $(D/\left(1-\frac{1}{n}\right)[\{0\}])$
to  $(D/\left(1-\frac{1}{m}\right)[\{0\}])$
iff $f'(z)\ne 0$ for $z\ne 0$ and $n \mult_0(f)=m$.

Examples of \'etale orbifolds morphisms are given in \S\ref{et-ex} below.

\section{Unfolding Orbicurves}

\begin{theorem}\label{thm-unfolding}
Let $(C/\Delta)$ be a smooth orbifold with $\dim(C)=1$.

Then there exists a finite \'etale
(in th sense of def.~\ref{def-etale})
orbifold morphism from a
curve $C'$ to $(C/\Delta)$, unless
$(C/\Delta)$ is isomorphic to 
\[
(\P_1/(1-\frac{1}{m})[\{\infty\}])
\]
or
\[
(\P_1/(1-\frac{1}{m})[\{\infty\}]+(1-\frac{1}{n})[\{0\}])
\]
with $m\ne n$.
\end{theorem}

As explained in \cite{Na} (Thm.~1.2.15), this follows from 
group-theoretical work
of Fox (\cite{F}), Bundgaard and Nielsen (\cite{BN}).

\subsection{Examples}\label{et-ex}
Consider the case $(2,2,2,2)$ (meaning that the support of $\Delta$ consists of $4$ distinct points with weights $1/2$ and multiplicity $2$ each).
For every four distinct points $p_i$ on $\P_1$ there exists
an elliptic curve $E$ with a $2\!:\!1$-ramified covering $\pi:E\to\P_1$
which is ramified precisely over the $p_i$. 
This covering is \'etale in the orbifold sense, and provides an unfolding of the given orbifold on $\P^1$.

Observe that  $\Aut(\P_1)$ acts triply transitively on $\P_1$, so that if the support of $\Delta$ consists of three points, these can be assumed to be $0,1,\infty$.

For the multiplicities $(2,4,4)$, $(2,3,6)$ and $(3,3,3)$ such
an unfolding can be obtained quite explicitly:

For the multiplicities $(2,4,4)$ we use the elliptic
curve $C$ defined by $y^2=x^3-x$ with ramified covering
$C\to\P_1$ given by the meromorphic function $x^2$.
Then above $0$ (resp. $1$, $\infty$), there are $1$ (resp. $2$; $1$) points with ramification multiplicities $4$ (resp. $2$; $4$), and no other ramification. This ramified cover is thus an unfolding of this $(2,4,4)$ orbifold on $\P^1$.

For the multiplicities $(2,3,6)$ we use the elliptic
curve $C$ defined by $y^2=x^3+1$ with ramified covering
$C\to\P_1$ given by the meromorphic function $y^2=x^3+1$.
Then above $0$ (resp. $1$, $\infty$), there are $3$ (resp. $2$; $1$) points with ramification multiplicities $2$ (resp. $3$; $6$), and no other ramification. This ramified cover is thus an unfolding of this $(2,3,6)$ orbifold on $\P^1$.

For the multiplicities $(3,3,3)$ we use the elliptic
curve $C$ defined by $y^2=x^3+1$ with ramified covering
$C\to\P_1$ given by the meromorphic function $y$.
Then above $-1$ (resp. $1$, $\infty$), there is one single point with ramification multiplicity $3$, and no other ramification. This ramified cover is thus an unfolding of this $(3,3,3)$ orbifold on $\P^1$.

\section{Fundamental group}

\begin{definition}\label{def-fd-gp}
Let $(X,\Delta)$ be an orbifold.

The {\em orbifold fundamental group} is the quotient
of $\pi_1(X\setminus|\Delta|)$ by the normal subgroup $N$
generated by all loops who can be realized as the image
of $t\mapsto \frac{1}{2}e^{2\pi i t}$ under some
classical orbifold morphism from the unit disk
$(D,\emptyset)$ to $(X,\Delta)$.
\end{definition}

\begin{lemma}
Assume that $X$ is smooth.
Then $N$ is generated by small loops around each connected
component of the smooth part of $|\Delta|$.
\end{lemma}
\begin{proof}
Let $H:[0,1]\times D\to X$ be a homotopy between $f:D\to X$
and a constant map with value $p\in X\setminus|\Delta|$.
Since $X$ is smooth we may assume by transversality arguments
that $H$ stays away from $\Sing(D)$.
This implies the statement.
\end{proof}

%

\begin{proposition}
A classical orbifold morphism induces a group homomorphism
between the orbifold fundamental groups.
\end{proposition}

\begin{remark}
This statement is very false for the (non-classical) orbifold
morphisms. For example, $z\mapsto z^d$ induces a non-classical
orbifold morphism from 
$(D/(1-\frac{1}{n})[\{0\}])$
to
$(D/(1-\frac{1}{m})[\{0\}])$
whenever $dn\ge m$ but there is no natural group homomorphism
from
$\pi_1(D/(1-\frac{1}{n})[\{0\}])=\Z/n\Z$ to
$\pi_1(D/(1-\frac{1}{m})[\{0\}])=\Z/m\Z$ unless $n$
divides $m$.
\end{remark}

\begin{proof}
Each element $\gamma\in\pi_1(X\setminus|\Delta|)$
can be represented by a loop inside
$X\setminus(|R_f + \Delta|)$. Let $\gamma_i$ ($I=1,2$) be such loops
homotopic to $\gamma$. Observe that $f(\gamma_i)\subset
X\setminus|\Delta'|$.
The $\gamma_i$
homotopic to each other in $X\setminus|\Delta|$.
For $x\in X\setminus|\Delta|$ we have $f(x)\not\in|\Delta'|$
unless $x\in R_f$. Hence the homotopy classes
of $f\circ\gamma_i$ differ only by an element of $N'$.

It follows that there is a group homomorphism between the orbifold
fundamental groups.
\end{proof}

\begin{proposition}\label{etale-lift}
Let $f:(X/\Delta)\to (X'/\Delta')$
be an \'etale orbifold morphism between smooth orbifold curves
 and let $g:(D/\emptyset)\to(X'/\Delta')$ be a classical orbifold
morphism.

Then there exists a classical orbifold morphism $\tilde g:
D\to (X/\Delta)$
such that $g=f\circ\tilde g$.
\end{proposition}

\begin{proof}
Local calculations verify that such lifts $\tilde g$ exist locally.
These local solutions then define a local system which is globally
trivial, because the disc is simply-connected. Hence there is a global lift \
$\tilde g$.
\end{proof}

\begin{remark}
Again this is very false for non-classical orbifold morphisms:
$h:z\mapsto z^n$ defines an \'etale orbifold morphism from
$D$ to $(D/(1-\frac{1}{n})[\{0\}])$, but for a given orbifold
morphism $g:D\to (D/(1-\frac{1}{n})[\{0\}])$ there exists
a lift $\tilde g$ only if $g$ is in fact a classical orbifold
morphism.
\end{remark}

\begin{proposition}\label{exist-et}
Let $(X/\Delta)$ be a smooth orbifold curve.
Let $\Gamma$ be a subgroup
of the orbifold fundamental group  $\pi_1(X/\Delta)$.

Then there exists an orbifold $(X'/\Delta')$ and an \'etale
orbifold map $f:(X'/\Delta')\to(X/\Delta)$ such that
$(\pi_1(f))(\pi_1(X'/\Delta'))=\Gamma$.
\end{proposition}

\begin{proof}
Recall that $\pi_1(X/\Delta)=\pi_1(X\setminus|D|)/N$ where $N$ is defined as
in def.~\ref{def-fd-gp}. 
Thus we obtain a subgroup $\Gamma_0\subset \pi_1(X\setminus|D|)$
such that $N\subset\Gamma_0$ and $\Gamma_0/N=\Gamma$. Let
$\rho:Y\to X\setminus|D|$ be the unramified covering of $X\setminus|D|$
associated to the subgroup $\Gamma_0\subset\pi_1(X\setminus|D|)$.
Consider now $p\in|D|$. We may embedd a small disc $D$ into $X$
such that $0$ is mapped to $p$ by the embedding map $i$. Then 
$\Lambda_p=i_*(\pi_1(D\setminus\{0\}))$ is a cyclic subgroup of
$\pi_1(X\setminus|D|)$ containing $N\cap \Lambda_p$ as subgroup
of finite index. Since $N\subset\Gamma_0$, we may deduce that
$\Gamma_0\cap\Lambda_p$ is of finite index in $\Lambda_p$.
It follows: $\rho^{-1}(i(d\setminus\{0\}$ decomposes into connected
component on each of which $\rho$ is isomorphic to $z\mapsto z^k$
where $k=[\Lambda_p:\Lambda_p\cap\Gamma_0]$.
Thus we can complete $Y$ over $p$ by adding one point to each
connected component of $\rho^{-1}(i(D\setminus\{0\}$.
The orbifold multiplicity for each of this added points has to be chosen
as $m/k$ where $m=[\Lambda_p:\Lambda_p\cap N]$ is the multiplicity for $p$.
Doing this for every point $p\in|D|$, we obtain an orbifold $(X'/\Delta')$
with an orbifold projection morphism $f:(X'/\Delta')\to(X/\Delta)$. 
\end{proof}

\begin{remark}\label{lift} Thus classical orbifold morphisms from the unit disc to an orbifold curve $(C/\Delta)$ can be lifted to unfoldings of $(C/\Delta)$, while their non classical versions cannot. For this reason the study of these classical maps reduces to the non orbifold case on any unfolding, while the study of the no classical version poses (seemingly) new problems. On the level of arithmetics, exactly the same situation appears: see \cite{D} for the classical orbifold version of Mordell's conjecture on curves, and \cite{C'} for its non-classical version (which is presently only a conjecture).
\end{remark}

However for the category of classical orbifold morphisms 
we obtain a Galois theory for coverings:

\begin{proposition}\label{etale-galois}
Let $\pi:(X/\Delta)$ be a smooth orbifold curve.

Then there is a natural one-to-one correspondance between
\begin{itemize}
\item
subgroups of $\Gamma$ of $\pi_1 (X/\Delta)$
\item
\'etale orbifold coverings $(X'/\Delta')\to(X/\Delta)$.
\end{itemize}
\end{proposition}
\begin{proof}
If $\Gamma$ is a subgroup of $\pi_1 (X/\Delta)$, the existence of
a corresponding \'etale covering follows from prop.~\ref{exist-et}.
Conversely let $\pi:(X'/\Delta')\to(X/\Delta)$ be an \'etale orbifold
cover. Let $N$ resp.~$N'$ be the subgroups of $\pi_1(X\setminus|D|)$
resp.~ $\pi_1(X'\setminus|D'|)$ as in def.~\ref{def-fd-gp}. 
Since $X\setminus|D|\to
X'\setminus |D'|$ is an unramifid covering, we obtain an embeding
of $\pi_1(X'\setminus|\Delta'|)$ into $\pi_1(X\setminus|\Delta|)$.
Due to prop.~\ref{etale-lift}, this embedding identifies $N$ with $N'$.
Hence the statement.
\end{proof}

As a consequence, for every smooth orbifold curve $(X/\Delta)$ there is
a smooth orbifold curve $(X'/\Delta')$ with $\pi_1(X'/\Delta')=\{e\}$
and a properly discontinuos action of $\Gamma=\pi_1(X/\Delta)$ on
$(X'/\Delta')$ such that $(X/\Delta)$ can be regarded as the quotient
of $(X'/\Delta')$ by this $\Gamma$-action.

\section{Uniformization}
\begin{proposition}\label{prop-uniformization}
A smooth one-dimensional orbifold $(X/\Delta)$ has trivial
fundamental group (as defined in def.~\ref{def-fd-gp}) if and only if it is
isomorphic to one of the following:
$\C$, $D$, $\P_1$ $(\P_1/(1-\frac{1}{n})\{\infty\})$
or $(\P_1/(1-\frac{1}{n})\{\infty\}+
(1-\frac{1}{m})\{0\})$ with $gcd(n,m)=1$.

For every smooth one-dimensional orbifold $(X/\Delta)$ there
exists a smooth one-dimensional orbifold $(\tilde X/\tilde\Delta)$
with trivial fundamental group and an \'etale orbifold morphism
$\pi:(\tilde X/\tilde\Delta)\to(X/\Delta)$.
\end{proposition}
\begin{proof}
The first statement follows from thm.~\ref{thm-unfolding} if $X$
is compact. If $X$ is not compact we note that $\pi_1(X/\Delta)=\{0\}$
implies that $X$ is simply-connected. Hence (in the non-compact case)
we have $X\simeq\C$ or $X\simeq D$. However for both $X=\C$ and $X=D$
it is immediate that $\pi_1(X/\Delta)\ne\{0\}$ unless $\Delta=0$.

The second statement follows from prop.~\ref{etale-galois}.
\end{proof}

\section{Hyperbolicity and Kobayashi pseudodistance}

We recall (and extend) from  \cite{C} the notion of orbifold Kobayashi pseudodistance by restricting to orbifold morphisms from the unit disc to $(X/ \Delta)$.

More precisely:

\begin{definition}
Let $(X/\Delta)$ be an orbifold
with $\Delta=\sum_i a_iH_i$.
Let $\Delta_1$ be the union of all $H_i$ with $a_i=1$ (ie: weight one, or equivalently multiplicity infinite).

The orbifold Kobayashi pseudodistance of the orbifold $(X/\Delta)$
is the largest pseudodistance on $(X\setminus|\Delta_1|)$ such that
every orbifold morphism from the unit disc $D$ to $(X/\Delta)$
is distance-decreasing with respect to the Poincar\'e distance on the
unit disc.

One defines similarly the classical orbifold Kobayashi pseudodistance on $(X/\Delta)$ by replacing the above set of orbifold morphisms from the disc to $(X/\Delta)$ by their classical versions.
\end{definition}

\begin{remark} Let $d_X$ (resp. $d_{(X/\Delta)}$; resp. $d^*_{(X/\Delta)})$ be the usual (resp. orbifold; resp. classical orbifold) Kobayashi pseudodistance. 
Then we have: 
\[
d_X\le d_{(X/\Delta)}\le d^*_{(X/\Delta)}
\le d_{X^\setminus|\Delta|}.
\]
It is clear that $d_X$ and  $d_{(X/\Delta)}$ are usually very different
as well as $d^*_{(X/\Delta)}$ and $d_{X\setminus|\Delta}$.
But we do not know a single example in which   $d_{(X/\Delta)}$
and $d^*_{(X/\Delta)})$ differ.
\end{remark}

The definition implies immediately that the (classical) orbifold
Kobayashi pseudodistance is distance-decreasing under (classical)
orbifold morphisms between orbifolds.

As in the case of the usual Kobayashi pseudodistance for
manifolds there is an equivalent definition using chains of disc:

{\em For $x,y\in X\setminus|\Delta_1|$ the (classical) Kobayashi pseudodistance
$d_{(X/\Delta)}$ is the infimum over $\sum_id_P(p_i,q_i)$ where
$d_P$ is the distance function on the unit disc $D$ induced
by the Poincar\'e metric and the infimum is taken over all finite families
$f_1,\ldots f_d$ of (classical) orbifold morphisms from $D$ to
$(X/\Delta)$ with $f_1(p_1)=x$, $f_d(q_d)=y$ and
$f_k(q_k)=f_{k+1}(p_{k+1})$.}

From this definition it is easily deduced that:

\[
d_{(X/\Delta)}: X\setminus|\Delta_1| \times X\setminus |\Delta_1| \to\R
\]
is continuous and that 
the set
\[
E_x=\{y\in X\setminus\Delta:d_{(X/\Delta)}(x,y)=0\}
\]
is connected 
for every $x\in X\setminus|\Delta_1|$.

\begin{definition}
An orbifold $(X/\Delta)$ is (classically) orbifold hyperbolic
if the (classical) orbifold Kobayashi pseudodistance is
a distance on $X\setminus\Delta_1$ where $\Delta_1$ is the union
of the components of $\Delta$ with multiplicity one.
\end{definition}

As a consequence of prop.~\ref{etale-lift} we obtain:

\begin{corollary}
Let $f:(X/\Delta)\to(X'/\Delta')$ be an \'etale orbifold
morphism between ondimensional orbifolds. Then $(X/\Delta)$ is classical  orbifold hyperbolic
if and only if $(X'/\Delta')$ has this property.
\end{corollary}

\section{Classical orbifold Kobayashi pseudodistance in dimension one}
\begin{proposition}
Let $(X/\Delta)$ be a one-dimensional smooth orbifold.

If there exists an \'etale orbifold morphism $\pi:D\to(X/\Delta)$,
then 
\[
d^*_{(X/\Delta)}(p,q)
=\inf_{x\in\pi^{-1}(p);y\in\pi^{-1}(q)}d_D(x,y)
\]

If there is no \'etale orbifold morphism $\pi:D\to(X/\Delta)$,
then $d^*_{(X/\Delta)}\equiv 0$.
\end{proposition}

\begin{proof}
Consequence of prop.~\ref{etale-lift} 
and prop.~\ref{prop-uniformization}.
\end{proof}

\begin{corollary}
Let $(X/\Delta)$ be a compact smooth one-dimensional orbifold.

Then $(X/\Delta)$ is classically hyperbolic iff $\deg(K_{(X/\Delta)})>0$.
\end{corollary}

\subsection{Examples}
We consider $X=D$, $\Delta=\left(1-\frac{1}{n}\right)[\{0\}]$.
Then $z\mapsto z^n$ yields an unfolding $D\to (X/\Delta)$ and
consequently the classical Kobayashi pseudodistance on $(X/\Delta)$
is the distance function induced by the ``push-forward'' of the
Poincar\'e metric on $D$ which is easily calculated as
\[
\frac{4dzd\bar z}{n^2|z|^{2-\frac{2}{n}}\left(1-|z|^{\frac{2}{n}}\right)^2}.
\]
Note that for $n\to\infty$ this converges
to
\[
\frac{4dzd\bar z}{|z|^2\left(\log|z^2|\right)^2}
\]
which is the push-forward of the Poincar\'e metric under the
universal covering map from $D$ to the punctured disc
$D^*=\{z\in\C:0<|z|<1\}$.
\section{An orbifold Brody theorem}
Brody's theorem (\cite{B}) is an important tool in the study
of hyperbolicity questions for complex spaces.
Here we will develop a version of this theorem for orbifolds.

As a first step we show:

\begin{proposition}\label{orbimaps-closed}
Let $f_n:(X/\Delta)\to(X'/\Delta')$ be a sequence of
orbifold morphisms.
Assume that $(f_n)$, regarded as a sequence of holomorphic maps from $X$ to $X'$
converge locally uniformly to a holomorphic map $f:X\to X'$.

Then either $f(X)\subset|\Delta'|$ or $f$ is an orbifold
morphism from $(X/\Delta)$ to $(X'/\Delta')$.

This statement, and its proof, hold both in the classical and non classical versions.
\end{proposition}

\begin{proof}
Assume $f(X)\not\subset|\Delta'|$.

Fix an orbifold morphism $g:D\to(X/\Delta)$.
By definition, $f_n\circ g$ are orbifold morphisms and we have
to show that $f\circ g$ is an orbifold morphism as well.

Let $D_i$ be an irreducible component of $\Delta$ with
multiplicity $\frac{m-1}{m}$.
Let $p\in D$ with $q=f(g(p))\in|D_i|$.
We have to show that $(f\circ g)^*D_i$ has multiplicity at
least $m$. In an
open neighbourhood $U$ of $q$ in $X$ the divisor $D_i$ has a
defining function $\rho$.
Let $W$ be a relatively
compact open neighbourhood of $p$ in $(f\circ g)^{-1}(U)$.
The set of all maps $F:X\to X'$ with $F(g(\bar W))\subset U$
is open for the the topology of locally uniform convergence.
Thus we have $f_n(g(W))\subset U$ for all sufficiently
large $n$. Now $\rho\circ f_n\circ g$ is a sequence of
holomorphic functions on $W$ converging to $\rho\circ f\circ g$.
Since we assumed that $f(X)$ is not contained in $|\Delta'|$,
$\rho\circ f\circ g$ does not vanish identically. Hence there
is a number $\epsilon>0$ such that $S_\epsilon(p)=
\{z\in \C:|z-p|=\epsilon\}$ is contained in $W$ and
$\rho\circ f \circ g$ has no zero in
$B_\epsilon(p)=\{z\in \C:|z-p|\le\epsilon\}$ except at $p$.

The theorem of Rouch\'e now implies that
for all sufficiently large $n$ the multiplicity $\mu$ of
$\rho\circ f\circ g$ at $p$ equals the sum of all multiplicities
of all zeroes in $B_\epsilon(r)$ of $\rho\circ f_n\circ g$.

Hence there is at least one zero  of $\rho\circ f_n\circ g$
in $B_\epsilon(r)$ for $n$ sufficiently large (since $f(p)\in|D_i|$).
Furthermore each such zero has multiplicity at least $m$, because
$f_n\circ g:D\to (X'/\Delta')$ are orbifold morphisms.

Therefore $\mu$ is at least $m$.
Since this argument may be applied to all components $D_i$
of $|\Delta|$ and all points $p\in D$ with $f\circ g(p)\in|D_i|$
for every orbifold morphism $g:D\to (X/\Delta)$,
we may conclude that $f$ is an orbifold morphism.
\end{proof}

\begin{remark}
As said, this works as well for both ``classical'' and ``non classical''
orbifold morphisms: In the last case we use the ordinary ordering
on $\N$ while in the first case we use the partial ordering of $\N$ by divisibility.
\end{remark}

\begin{proposition}\label{infini-nonhype}
Let $(X/\Delta)$ be an orbifold and let $\Delta_1$ be the union
of components of $\Delta$ with weight one (or equivalently, multiplicity $\infty$).
Assume that there are two distinct points $p,q\in X\setminus|\Delta_1|$
with orbifold Kobyashi pseudodistance zero.
Let $h$ be a hermitian metric on $X$ and let $d_h$  be the induced
distance function.

Then there exists a sequence of points $p_n\in X\setminus|\Delta_1|$
and orbifold morphisms $f_n:D\to (X/\Delta)$ such that
$f_n(0)=p_n$, $\lim p_n=p$ and $\lim||f_n'(0)||=+\infty$,
the latter calculated with respect to the Poincar\'e metric on $D$
and the hermitian metric $h$ on $X$.
\end{proposition}

\begin{proof}
If not, there exists a neighbourhood $W$ of $p$ and a constant $C>0$
such that $||f'(0)||\le C$ for all orbifold morphisms $f:D\to
(X/\Delta)$
with $f(0)\in W$. Let us assume that this is the case. Since $D$ is
homogeneous and the composition
$f\circ\phi$ is an orbifold morphism for every orbifold morphism
$f$ and every automorphism $\phi$ of $D$, this condition implies
that $||f'(z)||\le C$ for every orbifold morphism $f:D\to (X/\Delta)$
and every $z\in D$ with $f(z)\in W$. By shrinking $W$, we may assume
$q\not\in W$. Now for every $\epsilon>0$ there is a chain of
orbifold discs as in \S9 above with $\sum d_P(p_i,q_i)\le\epsilon$.
By taking geodesics in $D$ linking $p_i$ with $q_i$ and concatenating
their images we obtain a piecewise smooth path $\gamma:[0,1]\to X$
with $\gamma(0)=p$ and $\gamma(1)=q$.
Let $\alpha=\inf\{t:\gamma(t)\not\in W\}$.
Then
\[
\epsilon\ge d(p,\gamma(\alpha))\ge C d_h(p,\partial W)
\]
which leads to a contradiction since $d_h(p,\partial W)>0$.
\end{proof}

We recall the ``reparametrization lemma'' of Brody which may be
rephrased as follows:

\begin{proposition}\label{brody-re}
Let $X$ be a compact complex manifold and $f_n:D\to X$ a sequence
of holomorphic maps with $\limsup||f_n'(0)||=+\infty$.

Then there exists an increasing sequence of positive real numbers $r_n$
and a sequence of holomorphic maps $\alpha_n:D(r_n;0)\to D$
such that $\lim r_n=+\infty$ and such that a subsequence of
$f_n\circ\alpha_n$ converges locally uniformly to a holomorphic
map $f:\C\to X$ with
\[
\sup_{z\in\C}||f'(z)||=||f'(0)||>0.
\]
\end{proposition}

\begin{theorem}
Let $(X/\Delta)$ be a compact orbifold. Assume that the (classical) orbifold
Kobayashi pseudodistance on $X\setminus|\Delta_1|$ is not a distance.

Then there exists a non-constant holomorphic map $f:\C\to X$
which is either a (classical) orbifold morphism or fulfills
the property $f(\C)\subset|\Delta|$.

Furthermore
\[
\sup||f'(z)||=||f'(0)||>0.
\]
\end{theorem}

\begin{proof}
By prop.~\ref{infini-nonhype} there is a sequence of orbifold morphisms
$f_n:D\to (X/\Delta)$ such that $\lim||f'(0)||=+\infty$.
Due to ``Brody reparametrization'' (prop.~\ref{brody-re})
there are sequences $r_n\in\R^+$ and $\alpha_n:D(r_n;0)\to D$
such that $\lim r_n=+\infty$ and such that a subsequence
of $f_n\circ\alpha_n$ converges to a holomorphic map $f:\C\to X$
with $f'(0)\ne 0$. Now compositions of orbifold morphisms are
orbifold morphisms, hence $f_n\circ\alpha_n$ are orbifold morphisms.
Thus prop.~\ref{orbimaps-closed} implies that for all $r>0$
either $f|_{D_r}:D_r\to (X/\Delta)$ is an orbifold morphisms or
$f(D_r)\subset |\Delta|$.
As a consequence, either $f:\C\to (X/\Delta)$ is an orbifold
morphism or $f(\C)\subset|\Delta|$.
\end{proof}

\begin{corollary}\label{brody-curve}
Let $(X/\Delta)$ be a one-dimensional compact orbifold.

Then either $(X/\Delta)$ is orbifold hyperbolic or there exists
a non-constant orbifold morphism $f:\C\to (X/\Delta)$ with
bounded derivative.
\end{corollary}

\begin{proof}
Since $X$ is one-dimensional, $|\Delta|$ is discrete. As a consequence
$f(\C)$ can not be contained in $|\Delta|$ for a holomorphic map
$f:\C\to X$ with $f'(0)\ne 0$.
\end{proof}

\section{Nevanlinna theory}

We use the usual notations of Nevanlinna theory (see e.g. \cite{N}).
In particular, if $D$ is a divisor on a complex
space $X$ and $f:\C\to X$ is a holomorphic map, then
\[
N_f(r,D) = \int_1^r \deg(f^*D|_{D_t}) \frac{dt}{t}
\]
and
\[
N^1_f(r,D) = \int_1^r \deg((f^*D)_{red}|_{D_t}) \frac{dt}{t}
\]
where $D-t=\{z\in\C:|z|<t\}$.

If furthermore $\omega$ is a $(1,1)$-form on $X$ (e.g.~a K\"ahler form
or $c_1(L(D))$), then
\[
T_f(r,\omega)=\int_1^r \left ( \int_{D_t}f^*\omega \right ) \frac{dt}{t}.
\]

 \begin{proposition}\label{orbi-trunc}
Let $X$ be a compact complex manifold, $H$ an irreducible reduced
hypersurface, $n\in \N\cup\{+\infty\}$,
$\alpha=(1-1/n)$, $\Delta=\alpha H$ and
$f:\C\to (X/\Delta)$
an orbifold morphism.

Then
\[
T_f(r,c_1(H))-N_f^1(r,H)\ge \alpha T_f(r,c_1(H)).
\]
\end{proposition}

\begin{proof}
By the ``First Main Theorem'' of Nevanlinna theory, we have
\[
T_f(r,c_1(H))\ge N_f(r,H)\ge 0.
\]
Now $N^1_f(r,H)$ is the ``truncated counting function'' which ignores
multiplicities and $f^*H$ has multiplicity at least $n$
at every point of $f^{-1}|H|$.
Hence
\[
N_f(r,H)\ge n N_f^1(r,H)
\]
Together these two inequalities imply
\[
T_f(r,c_1(H))-N_f^1(r,H)\ge \left(1-\frac{1}{n}\right) T_f(r,c_1(H))
=\alpha T_f(r,c_1(H)).
\]
\end{proof}

\begin{definition}
We say that the ``S.M.T.\footnote{``S.M.T''=Second Main Theorem}
 with truncation level 1'' holds for
a holomorphic map $f$ from $\C$ to a compact complex manifold $X$
and a reduced effective divisor $D$ on $X$
if
\[
T_f(r,c_1(D+K))-N_f^1(r,D)\le \epsilon T_f(r,\omega)||_\epsilon
\]
for some positive $(1,1)$-form $\omega$ on $X$. (The notation $||_\epsilon$ means that the inequality holds for any $\epsilon >0$, for $r$ outside a subset of finite measure depending on $\epsilon$).
\end{definition}

By a classical result of Nevanlinna (\cite{N}), the
``S.M.T. with truncation level one''
holds for every non-constant holomorphic map to a one-dimensional
compact complex manifold $X$ and every reduced effective divisor $D$.

\begin{proposition}
Let $(X/\Delta)$ be a compact orbifold, and let $f:\C\to (X/\Delta)$
be an orbifold morphism such that  the ``S.M.T. with truncation level one''
holds for the underlying holomorphic map $f:\C\to X$ and the divisor
$H$ on $X$ which is obtained by replacing all mutiplicities by one.
Then
\[
T_f(r,c_1(\Delta+K_X))\le \epsilon T_f(r,\omega)||_\epsilon
\]
for every  positive $(1,1)$-form $\omega$ on $X$.
\end{proposition}

\begin{proof}
Let $\Delta=\sum\alpha_i H_i$ with $\alpha_i=1-\frac{1}{n_i}$.
Then $H=\sum H_i$.
Due to the S.M.T. we have:
\[
T_f(r,c_1(K_X))
+\sum_i \left( T_f(r,c_1(H_i))-N_f^1(r,H_i) \right)
\le \epsilon T_f(r,\omega)||_\epsilon
\]
By prop.~\ref{orbi-trunc}:
\[
T_f(r,c_1(H_i))-N_f^1(r,H_i) \ge \alpha_i T_f(r,c_1(H_i))
\]
Therefore:
\begin{align*}
T_f(r,c_1(\Delta+K_X))
&=T_f(r,c_1(K_X))+\sum_i \alpha_i T_f(r,c_1(H_i)) \\
& \le
T_f(r,c_1(K_X))
+\sum_i \left( T_f(r,c_1(H_i))-N_f^1(r,H_i) \right)\\
& \le
T_f(r,c_1(K_X))
+T_f(r,c_1(H))-N^1_f(r,c_1(H)) \\
&\le \epsilon T_f(r,\omega)||_\epsilon
\end{align*}
\end{proof}

\begin{corollary}\label{nevan-dgK}
Let $X$ be a compact smooth complex curve (i.e.~a compact Riemann surface)
of genus $g$ such that there exists  a non-constant
orbifold morphism  $f:\C\to(X/\Delta)$.

Then $\deg(\Delta+K_X)\le 0$, i.e. $\deg(\Delta)\le 2-2g$.
\end{corollary}

\begin{proof}
For curves, the ``S.M.T. with truncation level one''
has already been established by Nevanlinna (\cite{N}).
It follows that $\deg(\Delta+K_X)\le 0$ whenever there exists
a non-constant orbifold morphism. But $\deg(\Delta+K_X)\le 0$
is equivalent to $\deg(\Delta)\le -\deg K_X=2-2g$.
\end{proof}

\section{Hyperbolicity of orbicurves}

We characterize completely under which condition an orbifold
of dimension one is orbifold hyperbolic.

\begin{theorem}
Let $(X/\Delta)$ be a smooth orbifold curve.

Then $(X/\Delta)$ is orbifold hyperbolic if and only if it is
classically orbifold hyperbolic.

If $X$ can be compactified
to a smooth compact curve $\bar X$ by adding finitely many points
and in addition the support $|\Delta|$ is finite, then the
orbifold hyperbolicity of $(X/\Delta)$ is equivalent to
$\deg(K_{\bar X}+\Delta)+\#(\bar X\setminus X)>0$.

Otherwise (if there is no such compactification or the support
$|\Delta|$ is infinite) the orbifold $(X/\Delta)$ is
orbifold hyperbolic.
\end{theorem}

\begin{proof}
We recall that a Riemann surface $X$ is hyperbolic unless
it is an elliptic curve, $\P_1$, $\C$ or $\C^*$. In particular,
if $X$ can not be compactified by adding finitely points,
it must be hyperbolic and as a consequence $(X/\Delta)$ is
orbifold hyperbolic and classically orbifold hyperbolic.

Now assume that $|\Delta|$ is finite and $X$ can be compactified
by adding finitely many points. By adding these points to $\Delta$
(with weigth $1$)
we may assume that $X$ is already compact.
If $(X/\Delta)$ is not hyperbolic, there is a orbifold morphism
from $\C$ to $(X/\Delta)$ due to cor.~\ref{brody-curve}. 
Using Nevanlinna theory (see cor.~\ref{nevan-dgK}), 
this implies
$\deg(K_X+\Delta)\le 0$.
On the other hand, if $\deg(K_X+\Delta)\le 0$, there are two
possibilities: Either $X$ is an elliptic curve and $\Delta$ is empty
or $X\simeq\P_1$. Evidently elliptic curves are not hyperbolic.
Thus it remains to discuss the case $X=\P_1$. If $|\Delta|$ contains
at most two points, $\C^*$ embedds into $(X/\Delta)$ which therefore
can not be hyperbolic. Finally, if $|\Delta|$ contains at least
three points, due to thm.~\ref{thm-unfolding} 
there is an \'etale orbifold morphism from
a compact curve $C$ to $(X/\Delta)$. Now $\deg(K_X+\Delta)\le 0$
implies $\deg(K_C)\le 0$ and thereby implies that is either $\P_1$
or an elliptic curve. In both cases the projection map from $C$ to
$(X/\Delta)$ shows that the latter is not hyperbolic.

We still have to discuss the case where $X$ can be compactified by
adding finitely many points, but $|\Delta|$ is infinite.
Because $|\Delta|$ is infinite and the multiplicity at each
point is at least $\frac{1}{2}$,
we can find a finite $\Qp$-Weil divisor
$\Delta'$
by taking finitely many components of $\Delta$ with the same multiplicities
in such a way that $\deg(\Delta')$ is as large as desired.
Therefore there is a finite $\Qp$-Weil divisor $\Delta'$ on $X$
such that
\begin{enumerate}
\item the identity map of $X$ gives a classical orbifold morphism
from $(X/\Delta)$ to $(X/\Delta')$
\item $\deg(K_{\bar X}+\Delta')>0$.
\end{enumerate}
It follows that $(X/\Delta)$ is classically orbifold hyperbolic and therefore
orbifold hyperbolic.
\end{proof}

\begin{corollary}
Let $(X/\Delta)$ be a one-dimensional smooth compact orbifold.

Then $(X/\Delta)$ is not orbifold hyperbolic if and 
only if one of the following conditions
hold:
\begin{enumerate}
\item
$X$ is an elliptic curve and $\Delta$ is empty.
\item
$X\simeq\P_1$ and $|\Delta|$ contains at most two points.
\item
$X\simeq\P_1$ and there are
numbers: $p\le q\le r\in\N\cup\{\infty\}\setminus\{1\}$
such that $(X/\Delta)$ is isomorphic to
\[
\left(\P_1/(1-\frac{1}{p})\{0\}+(1-\frac{1}{q})\{1\}+(1-\frac{1}{r})\{\infty\}
\right)
\]
and $1/p+1/q+1/r\ge 1$. (There are exactly $5$ possibilities for $(p,q,r)$: $(2,3,4); (2,3,5);(2,3,6);(2,4,4);(3,3,3)$).
\item
There is a point $\lambda\in\C\setminus\{0,1\}$ such that
 $(X/\Delta)$ is isomorphic to
\[
\left(\P_1/(1-\frac{1}{2})\{0\}+
(1-\frac{1}{2})\{1\}+(1-\frac{1}{2})\{\infty\}+(1-\frac{1}{2})\{\lambda\}
\right)
\]
\end{enumerate}
\end{corollary}
\begin{proof}
A case-by-case check verifies that these are exactly the orbifold
curves for which $\deg(\Delta+K_X)\le 0$.
\end{proof}

\begin{remark} \label{classvsnclass}
The situation is actually much less understood as may appear.

If $(X/\Delta)$ is an hyperbolic orbicurve, we do not know whether or not the classical and non classical Kobayashi pseudodistances coincide,
not even in the most simple case where 
$(X/\Delta)=(D/\left(1-\frac{1}{n}\right)[\{0\}])$ 
with $n\in\N\setminus\{1\}$.

In higher dimensions, is it still true that classical hyperbolicity coincides with (non classical) hyperbolicity ?

Is there any concrete example where one can calculate the 
(non-classical) orbifold
Kobayashi pseudodistances (if these are not degenerate)?

What about the arithmetic counterpart? For the ``classical'' variant
there is the work of Darmon (\cite{D}), but nothing seems to be
known about the ``non-classical'' variant.
\end{remark}

\end{document}